\documentclass[draftclsnofoot,onecolumn, 12pt]{IEEEtran}
\usepackage{amsmath}
\usepackage{graphicx}
\usepackage{amssymb}
\usepackage{bm}
\usepackage{cite}
\newtheorem{theorem}{Theorem}

\newtheorem{cor}{Corollary}

\begin{document}

\title{\huge A New Type Of Upper And Lower Bounds On Right-Tail Probabilities Of Continuous Random Variables}
  
\author{Nikola Zlatanov  
\thanks{N. Zlatanov is with  Innopolis University, Innopolis, 420500, Russia.   E-mail: n.zlatanov@innopolis.ru.} 
  }
\maketitle

\begin{abstract}
In this paper, I present a completely new type of upper and lower bounds on the right-tail probabilities of continuous random variables with unbounded support and with semi-bounded support from the left. The presented upper and lower right-tail bounds  depend only on the probability density function  (PDF), its first derivative, and two parameters that are used for tightening the bounds. These tail bounds hold under certain conditions that depend on the PDF, its first and second derivatives, and the two parameters. The new tail bounds are shown to be  tight for a wide range of continuous random variables via  numerical examples.
\end{abstract}

 \section{Introduction}
 
The most well known and the most used method  for bounding tail probabilities is based on  variations of Markov's inequality \cite{markov_book}.
Markov's inequality   relates the tail probability of a non-negative random variable (RV) to its mean. The Bienaymé-Chebyshev's inequality  \cite{bienayme1853considerations, chebyshev1867valeurs},   relates the tail probability of a  RV to its mean and variance, and this inequality can be obtained by Markov's inequality. Other notable  bounds on the tail probabilities that  are based on   Markov's inequality are the Chernoff-Cramér  bound \cite{chernov} and   Hoeffding's inequality \cite{hoeff}, among the most famous. 

Additional tail bounding methods include martingale methods  \cite{mcdiarmid_1989}, information-theoretic methods \cite{ahlswede1976bounds, 10.5555/1146355},  the entropy method  based on logarithmic Sobolev inequalities  \cite{ledoux_1997},  Talagrand’s
induction method \cite{talagrand1995concentration}, etc.
For an  overview of tail bounding  methods, please refer to \cite{con_in_book}. 

In this paper, I present a completely new type of upper and lower right-tail bounds of continuous RVs  with unbounded support and with semi-bounded support from the left, under the assumption that the RVs satisfy certain conditions. The presented upper and lower right-tail bounds  depend only on the probability density function  (PDF), its first derivative, and two parameters that are used for tightening the bounds. These tail bounds hold under certain conditions that depend on the PDF, its first and second derivatives, and the two parameters.  By evaluating these conditions, one is able to establish whether a corresponding bound holds or not. Moreover, the two parameters can be used for   reshaping these conditions such that they are met, and thereby make the corresponding bounds to hold. In the paper, I also analyse and discuss  the tightness of the bounds by using the convergence rate between the upper and the lower bounds. In addition, I propose methods for optimizing the two parameters such that the tightness of the bounds is maximized. Finally, in the numerical section, I numerically evaluate the presented right-tail bounds for some well known RVs  and illustrate the tightness of the bounds.

 \section{The Type of Continuous Random Variables}

 Let $X$ be a continuous RV. Let $F_X(x)$ and $f_X(x)$ be the cumulative distribution function (CDF) and  the PDF of $X$, given by
 \begin{align}
 F_X(x)&={\rm Pr}\{X\leq  x\}\label{eq_1},\\
  f_X(x)&=\frac{d  F_X(x)}{d  x} ,\label{eq_2}
 \end{align}
 where  ${\rm Pr}\{A\}$ denotes the probability of an event $A$.

  The probabilities ${\rm Pr}\{X\leq  x\}$ and ${\rm Pr}\{X\geq  x\}$ are known as the left-tail and  right-tail probabilities of $X$, respectively.

Let  $f'_X(x)$ and $f''_X(x)$ be the first and second derivative of the PDF, given by
  \begin{align}
  f'_X(x)&=\frac{d  f_X(x)}{d  x}  \label{eq_3}
 \end{align}
and
  \begin{align}
  f''_X(x)&=\frac{d ^2 f_X(x)}{d  x^2} , \label{eq_3a}
 \end{align} 
 respectively.
 
In the following, when I say the support of the RV $X$, I mean the support of its PDF $f_X(x)$. Hence, support of $X$ and support of  $f_X(x)$ will be used interchangeably. Moreover, I will use $[(x_0, \infty)$ to mean    $(x_0, \infty)$ or $[x_0, \infty)$.
 
The bounds in this paper are derived based on the following assumptions:
\begin{enumerate}
\item The PDF of the RV $X$, $f_X(x)$, has  support on  $([x_0,\infty)$, where $-\infty\leq  x_0<\infty$.
\item $f_X(x)$ is a continuous function of $x$ on the entire  support of $X$.
\item $f'_X(x)$ exists and  is a continuous function of $x$ on the entire   support of $X$.
\item $f''_X(x)$ exists and  is a continuous function of $x$ on the entire   support of $X$.
\vspace{3mm}
\item $\lim\limits_{x\to\infty} \dfrac{f^2_X(x)}{f'_X(x)}=0$.
\vspace{3mm}
\item $\lim\limits_{x\to x_0} x f_X(x)=0.$
\vspace{3mm}
\end{enumerate} 

Note that condition 1) implies the following
\begin{align}
\lim\limits_{x\to \infty} f_X(x)&=0,\label{eq_4a} \\
\lim\limits_{x\to \infty} x f_X(x)&=0.\label{eq_4a-1}
\end{align}

 \section{Upper Bounds}
In this section, I provide upper bounds on the right-tail probabilities of   RVs satisfying assumptions 1)-6). These upper bounds  hold under specified conditions. The first bound is given in the following theorem. 
 
 \begin{theorem}\label{thm_1}
Let $X$ be a RV that has support on $([0,\infty)$ and satisfies assumptions 2) to 6). Then, an upper bound on the right-tail probability of $X$ is given by
\begin{align} \label{eq_4}
{\rm Pr}\{X\geq  x\}\leq - \frac{x^a f_X^2(x)}{f_X(x)+x^a f'_X(x)},
\end{align}
when $a$,  $x$, $f_X(x)$, $f'_X(x)$, and $f''_X(x)$ satisfy the following inequalities 
\begin{align} 
&a>0,\label{eq_a}\\
&f_X(x)+x^a f'_X(x)<0,\label{eq_5}\\
&(x-a x^a) f^2_X(x) 
- x^{2a+1} \big(f'_X(x)\big)^2
+x^{a+1}f_X(x)\big(f'_X(x)+x^a f''_X(x)\big) \leq 0.\label{eq_5a}
\end{align}
 \end{theorem}

\begin{IEEEproof}
The proof is provided in Appendix~\ref{app_1}.
\end{IEEEproof}

In the following theorem, I provide a corresponding upper bound when the RV $X$ has  support on $([x_0,\infty)$, for any $-\infty <x_0<\infty$.

\begin{theorem}\label{thm_2}
Let $X$ be a RV that has support on $([x_0,\infty)$, for $-\infty <x_0<\infty$, and satisfies assumptions 2) to 6). Then, an  upper bound on its right-tail probability is given by
\begin{align} \label{eq_6a}
{\rm Pr}\{X\geq  x\} \leq - \frac{(x-x_0)^a f_X^2(x)}{f_X(x)+(x-x_0)^a f'_X(x)},
\end{align}
when $a$, $x_0$, $x$, $f_X(x)$, $f'_X(x)$, and $f''_X(x)$ satisfy the following inequalities 
\begin{align} 
&a>0,\label{eq_a-1}\\
 &f_X(x)+(x-x_0)^a f'_X(x)<0,\label{eq_6b}\\
&(x-x_0-a (x-x_0)^a) f^2_X(x) 
- (x-x_0)^{2a+1} \big(f'_X(x)\big)^2\nonumber\\
&\qquad+(x-x_0)^{a+1}f_X(x)\big(f'_X(x)+(x-x_0)^a f''_X(x)\big) \leq 0.\label{eq_6c}
\end{align}
\end{theorem}

\begin{IEEEproof}
The proof is provided in Appendix~\ref{app_2}.
\end{IEEEproof}

Note that Theorem~\ref{thm_2} is a generalization of Theorem~\ref{thm_1} and includes the result of Theorem~\ref{thm_1} when $x_0$ in  Theorem~\ref{thm_2} is set to $x_0=0$.

The parameter $a$ in the bound and in the conditions of Theorem~\ref{thm_2} is a parameter that can be optimized in order to make the bound in \eqref{eq_6a} tighter and/or to reshape  the  conditions  in \eqref{eq_6b} and \eqref{eq_6c} such that they are met and thereby the bound in \eqref{eq_6a} holds. However, setting the value of $a$ to $a=1$ or to $a\to\infty$ in Theorem~\ref{thm_2} leads to much simpler bounds, which is made precise in the following two corollaries.

 \begin{cor}\label{cor_1}
Let $X$ be a RV that has support on $([x_0,\infty)$, for $-\infty <x_0<\infty$, and satisfies assumptions 2) to 6). Then, an  upper bound on its right-tail probability is given by
\begin{align} \label{eq_6a-cor_2}
{\rm Pr}\{X\geq  x\} \leq - \frac{(x-x_0) f_X^2(x)}{f_X(x)+(x-x_0) f'_X(x)},
\end{align}
when $x_0$, $x$, $f_X(x)$, $f'_X(x)$, and $f''_X(x)$ satisfy the following inequalities 
\begin{align} 
&f_X(x)+(x-x_0) f'_X(x)<0,\label{eq_6b-cor_2}\\
&
f'_X(x)+(x-x_0) f''_X(x) \leq \frac{(x-x_0)\big(f'_X(x)\big)^2}{f_X(x)}.\label{eq_6c-cor_2}
\end{align}
\end{cor}

\begin{IEEEproof}
Setting $a=1$ in Theorem~\ref{thm_2} leads directly to Corollary~\ref{cor_1}.
\end{IEEEproof}

Now, setting $a\to\infty$ in \eqref{eq_6a}, would lead to an even simpler and tighter bound than that in Corollary~\ref{cor_1},  if the corresponding conditions hold. However, note that  this simplified and tighter bound can also be obtained if $a$ is fixed in \eqref{eq_6a} and  I let $x_0\to-\infty$. Hence, this simplified and tighter upper bound is much more general and it holds for any $-\infty \leq x_0<\infty$. I specify this bound in the following corollary.

\begin{cor}\label{cor_2}
Let $X$ be a RV that has support on $([x_0,\infty)$, for $-\infty \leq x_0< \infty$, i.e., including the support $(-\infty,\infty)$, and satisfies assumptions 2) to 6). Then, an  upper bound on its right-tail probability is given by
\begin{align} \label{eq_4-1}
{\rm Pr}\{X\geq  x\}\leq - \frac{f_X^2(x)}{ f'_X(x)},
\end{align}
when $f_X(x)$, $f'_X(x)$, and $f''_X(x)$ satisfy the following inequalities
\begin{align} 
f'(x)&<0,\label{eq_5a-1-0}\\
\dfrac{f_X(x) f''(x)}{\big(f'(x)\big)^2}&\leq 1.\label{eq_5a-1}
\end{align}
 \end{cor}

\begin{IEEEproof}
Two proofs are possible. Fix $x_0$ and let $a\to\infty$ in Theorem~\ref{thm_2}. Then \eqref{eq_6a} becomes \eqref{eq_4-1}, condition \eqref{eq_6b} is always satisfied if \eqref{eq_5a-1-0} holds, and condition \eqref{eq_6c} becomes condition \eqref{eq_5a-1}. The same result is obtained if one fixes $a$ to any number larger than one and lets $x_0\to-\infty$ in  Theorem~\ref{thm_2}. This concludes the proof.
\end{IEEEproof}

Notably, the bound  in Corollary~\ref{cor_2} has  the simplest form so far. But there is more to it than just the simplicity. Specifically,  the bound in Corollary~\ref{cor_2} is tighter than or equal to the  bound in Theorem~\ref{thm_2}, and thereby also to the one in Corollary~\ref{cor_1}, when the corresponding conditions in Corollary~\ref{cor_2} hold. This is because the bound in Theorem~\ref{thm_2} is a decreasing function of $a$ when $f'_X(x)<0$. Hence, the tightest form of the bound in Theorem~\ref{thm_2} is the bound in Corollary~\ref{cor_2}, when the corresponding conditions hold. Moreover, this is the only upper bound, among the given, that holds for RVs with unbounded support on both sides of the real line. Hence, for RVs with unbounded support on both sides of the real line, it is not possible to improve the tightness of the bound  in Corollary~\ref{cor_2} by optimizing $a$, simply because this bound has already been optimized with respect to $a$ (by setting $a\to\infty$).

One conjecture that I would like to make here is the following. For all RVs $X$ that satisfy assumptions 2) to 6) and have  support on $(-\infty,\infty)$, the conditions in \eqref{eq_5a-1-0} and \eqref{eq_5a-1} always hold  $\forall x\geq x_0$, where $x_0<\infty$.

 However, not always condition \eqref{eq_5a-1} in Corollary~\ref{cor_2}  holds for  RVs that have   support on $([x_0,\infty)$, for $-\infty<x_0<\infty$. Specifically, 
 there are RVs with  support on $([x_0,\infty)$, where $-\infty<x_0<\infty$, for which condition \eqref{eq_5a-1} in Corollary~\ref{cor_2} does not hold and yet conditions \eqref{eq_6b-cor_2} and \eqref{eq_6c-cor_2} jointly hold, and thereby the bound in Corollary~\ref{cor_1}, and more generally the bound in Theorem~\ref{thm_2}, hold. The opposite is never true.
 
Another observation that I would like to make in the form of a conjecture is the following. Assume that $X$ has  support on $([x_0,\infty)$, for $-\infty<x_0<\infty$, and that condition \eqref{eq_5a-1} in Corollary~\ref{cor_2} does not hold, and that conditions \eqref{eq_6b-cor_2} and \eqref{eq_6c-cor_2} in Corollary~\ref{cor_1} jointly hold. Then, the optimal value of the parameter $a$ that minimizes the gap between the tail and the upper bound in Theorem~\ref{thm_2} as $x\to\infty$ is $a=1$.
  
 I make this conjecture, since I was not able to find a distribution for which it does not hold.  If this conjecture is true, then the consequence would be that the optimal value of the parameter $a$ in Theorem~\ref{thm_2} is either $a\to\infty$ or $a=1$, when $x\to\infty$.

\section{Lower Bounds}

I am also able to provide corresponding lower bounds on the right-tail probability of $X$ under certain conditions. These bounds are introduced and discussed in this section.

\begin{theorem}\label{thm_3}
Let $X$ be a RV that has support on $([0,\infty)$ and satisfies assumptions 2) to 6). Then, a  lower bound on its right-tail probability is given by
\begin{align} \label{eq_6a-nn1}
{\rm Pr}\{X\geq  x\}\geq - \frac{x^b}{1+x^b}\frac{x^a f_X^2(x)}{f_X(x)+x^a f'_X(x)},
\end{align}
when $a$, $b$, $x$, $f_X(x)$, $f'_X(x)$, and $f''_X(x)$ satisfy the following inequalities
\begin{align} 
&a>0,\label{eq_a-3}\\
&b>0,\label{eq_b-1}\\
&f_X(x)+x^a f'_X(x)<0,\label{eq_6c-nn3-n1}\\
&\Big(x\big(1+x^b\big)^2
- x^{a+b}(a+b+a x^b) \Big)f^2_X(x)
-x ^{2a+1}(x^{2b}-1) \big(f'_X(x)\big)^2\nonumber\\
&+ x^a f_X(x)\Big[ x(1+x^b)(2+x^b)f'_X(x)
+x^{a+b} \Big(-b f'_X(x)+x(1+x^b)f''_X(x)\Big)\Big]\geq 0.\label{eq_6c-nn3}
\end{align}
\end{theorem}

\begin{IEEEproof}
The proof is provided in Appendix~\ref{app_3}.
\end{IEEEproof}

For a RV $X$ with support on  $([x_0,\infty)$, the following theorem is applicable.

\begin{theorem}\label{thm_4}
Let $X$ be a RV that has support on $([x_0,\infty)$, for $-\infty < x_0<\infty$, and satisfies assumptions 2) to 6). Then, a  lower bound on its right-tail probability is given by
\begin{align} \label{eq_6a-nn1sd}
{\rm Pr}\{X\geq  x\}\geq - \frac{(x-x_0)^b}{1+(x-x_0)^b}\frac{(x-x_0)^a f_X^2(x)}{f_X(x)+(x-x_0)^a f'_X(x)},
\end{align}
when $a$, $b$, $x_0$, $x$, $f_X(x)$, $f'_X(x)$, and $f''_X(x)$ satisfy the following inequalities 
\begin{align} 
&a>0,\label{eq_a-3-1}\\
&b>0,\label{eq_b-1-1}\\
&f_X(x)+(x-x_0)^a f'_X(x)<0,\label{eq_6c-nn3-n1-1}\\
&\Big(y\big(1+y^b\big)^2
- y^{a+b}(a+b+a y^b)\Big) f^2_X(x)
-y^{2a+1}(y^{2b}-1) \big(f'_X(x)\big)^2\nonumber\\
&+ y^a f_X(x)\Big[ y(1+y^b)(2+y^b)f'_X(x)
+y^{a+b} \big(-b f'_X(x)+y(1+y^b)f''_X(x)\big)\Big]\geq 0,\label{eq_6c-nn3-1}
\end{align}
where 
\begin{align}\label{eq_y_x}
y=x-x_0.
\end{align}

\end{theorem}

\begin{IEEEproof}
The proof is straightforward using the approach shown in the proof of Theorem~\ref{thm_2}, and is therefore omitted.
\end{IEEEproof}

Note that now two parameters, $a$ and $b$, can be optimized to make the bound in Theorem~\ref{thm_4} tighter and/or make the corresponding conditions to be met and thereby the bound to hold. If one chooses not to optimize with respect to the parameter $a$ by setting   $a=1$, then the following bound is applicable.

\begin{cor}\label{cor_3}
Let $X$ be a RV that has support on $([x_0,\infty)$, for $-\infty < x_0<\infty$, and satisfies assumptions 2) to 6). Then, a  lower bound on its right-tail probability is given by
\begin{align} \label{eq_6a-nn1sd-cor_4}
{\rm Pr}\{X\geq  x\}\geq - \frac{(x-x_0)^b}{1+(x-x_0)^b}\frac{(x-x_0) f_X^2(x)}{f_X(x)+(x-x_0) f'_X(x)},
\end{align}
 when  $b$, $x_0$, $x$, $f_X(x)$, $f'_X(x)$, and $f''_X(x)$ satisfy the following inequalities 
\begin{align} 
&b>0,\label{eq_b-1-1-cor_4}\\
&f_X(x)+(x-x_0) f'_X(x)<0,\label{eq_6c-nn3-n1-1-cor_4}\\
&\Big(y\big(1+y^b\big)^2
- y^{1+b}(1+b+ y^b) \Big) f^2_X(x)
-y^{3}(y^{2b}-1) \big(f'_X(x)\big)^2\nonumber\\
&+ y f_X(x)\Big[ y(1+y^b)(2+y^b)f'_X(x)
+y^{1+b} \big(-b f'_X(x)+y(1+y^b)f''_X(x)\big)\Big]\geq 0,\label{eq_6c-nn3-1-cor_4}
\end{align}
where $y$ is given by \eqref{eq_y_x}. 
\end{cor}

\begin{IEEEproof}
Setting $a=1$ in Theorem~\ref{thm_4} leads directly to Corollary~\ref{cor_3}.
\end{IEEEproof}

Now setting $a\to\infty$ in Theorem~\ref{thm_4} would not lead to a lower bound  similar to the upper bound in Corollary~\ref{cor_2} that holds for RVs with support on $(-\infty, \infty)$. Instead, it will only lead to a lower bound  for RVs with support on $([x_0, \infty)$, where   $-\infty<x_0<\infty$. This is due to the term
$(x-x_0)^b/(1+(x-x_0)^b)$  in \eqref{eq_6a-nn1sd}. Nevertheless, setting $a\to\infty$ in Theorem~\ref{thm_4} would still lead to the following useful corollary.

\begin{cor}\label{cor_4}
Let $X$ be a RV that has support on $([x_0,\infty)$, for $-\infty <x_0< \infty$, and satisfies assumptions 2) to 6).  Then, a lower bound on its right-tail probability is given by
\begin{align} \label{eq_4-1dfsd}
{\rm Pr}\{X\geq  x\}\geq - \frac{(x-x_0)^b}{1+(x-x_0)^b} \frac{f_X^2(x)}{ f'_X(x)},
\end{align}
 when $b$, $x_0$, $x$, $f_X(x)$, $f'_X(x)$, and $f''_X(x)$ satisfy the following inequalities 
\begin{align} 
&b>0,\label{eq_b-1-1-1}\\
&f'_X(x)<0,\label{eq_6c-nn3-n1-1-1}\\
& 1-(x-x_0)^{2b}+
\frac{f_X(x)}{\big(f'_X(x)\big)^2}(x-x_0)^{b-1}
\Big(
-b f'_X(x)+\big(1+(x-x_0)^b\big) (x-x_0) f''_X(x)\Big)
\geq 0 .\label{eq_6c-nn3-1-1}
\end{align}
\end{cor}

\begin{IEEEproof}
 Fix $x_0$ and let $a\to\infty$ in Theorem~\ref{thm_4}. Then \eqref{eq_6a-nn1sd} becomes \eqref{eq_4-1dfsd}, condition \eqref{eq_6c-nn3-n1-1} is always satisfied if \eqref{eq_6c-nn3-n1-1-1} holds, and condition \eqref{eq_6c-nn3-1} becomes condition \eqref{eq_6c-nn3-1-1}. This concludes the proof.
\end{IEEEproof}

I am now left with
establishing a lower bound for RVs with support on $(-\infty,\infty)$. However, such a lower bound requires taking the mean into an account. This is made precise in the following theorem.

\begin{theorem}\label{thm_5}
Let $X$ be a RV that has support on $(-\infty,\infty)$,  and satisfies assumptions 2) to 6). Moreover, let the mean of $X$ be denoted by $\mu$, which is given by
$$\mu=\int_{-\infty}^\infty x f_X(x) dx.$$
 Then, a lower bound on the right-tail probability of $X$ is given by
\begin{align} \label{eq_4-1dfsd-1}
{\rm Pr}\{X\geq  x\}\geq - \frac{(x-\mu)^b}{1+(x-\mu)^b} \frac{f_X^2(x)}{ f'_X(x)},
\end{align}
when $b$, $\mu$, $x$, $f_X(x)$, $f'_X(x)$, and $f''_X(x)$ satisfy the following inequalities
\begin{align} 
&b>0,\label{eq_b-1-1-1-thm_7}\\
& x>\mu\\
&f'_X(x)<0\label{eq_6c-nn3-n1-1-1-thm_7}\\
& 1-(x-\mu)^{2b}+
\frac{f_X(x)}{\big(f'_X(x)\big)^2}(x-\mu)^{b-1}
\Big(
-b f'_X(x)+\big(1+(x-\mu)^b\big)(x-\mu) f''_X(x)\Big)
\geq 0 .\label{eq_6c-nn3-1-1-thm_7}
\end{align}
\end{theorem}

\begin{IEEEproof}
The proof is provided in Appendix~\ref{app_4}.
\end{IEEEproof}

In the following, I write the upper and lower tail bounds   in a different form, a form which provides an intuitive interpretation.

If the conditions specified in Corollary~\ref{cor_2} and Theorem~\ref{thm_5} are jointly satisfied for a RV $X$ with support on $(-\infty,\infty)$, then the following holds
\begin{align} \label{eq_5a3sdsd-1}
\dfrac{d  }{d  x}\left[ \dfrac{1}{f_X(x)}\right]  \leq \dfrac{1}{{\rm Pr}\{X\geq  x\}}\leq \left(1+\frac{1}{(x-\mu)^b}\right) \dfrac{d  }{d  x}\left[ \dfrac{1}{f_X(x)}\right],
\end{align}
which means that the multiplicative inverse of the right-tail probability, ${\rm Pr}\{X\geq  x\}$,  of the RV $X$ is   bounded  from above and bellow by  the derivative of the multiplicative inverse of its  PDF, $f_X(x)$. 

From \eqref{eq_5a3sdsd-1}, it is easy to see  that the lower and upper bounds converge to each other with rate $1/(x-\mu)^b$, for $x>\mu$.

Next, if the conditions specified in Corollary~\ref{cor_2} and Corollary~\ref{cor_4} are jointly satisfied for a RV $X$ with support on $([x_0,\infty)$,  for $-\infty<x_0<\infty$, then the following holds
\begin{align} \label{eq_5a3sdsd-2}
\dfrac{d  }{d  x}\left[ \dfrac{1}{f_X(x)}\right]  \leq \dfrac{1}{{\rm Pr}\{X\geq  x\}}\leq \left(1+\frac{1}{(x-x_0)^b}\right) \dfrac{d  }{d  x}\left[ \dfrac{1}{f_X(x)}\right],
\end{align}
which again means that the multiplicative inverse of the right-tail probability, ${\rm Pr}\{X\geq  x\}$,  of the RV $X$ is  bounded   from above and bellow by  the derivative of the multiplicative inverse of its  PDF, $f_X(x)$. 

From \eqref{eq_5a3sdsd-2}, it is easy to see  that the lower and upper bounds converge to each other with rate $1/(x-x_0)^b$, for $x>x_0$.

Finally, if the conditions specified in Corollary~\ref{cor_1} and Corollary~\ref{cor_3} are jointly satisfied for a RV $X$ with support on $([0,\infty)$, then the following holds
\begin{align} \label{eq_5a3sdsd-3}
\dfrac{d  }{d  x}\left[ \dfrac{1}{x f_X(x)}\right]  \leq \dfrac{1}{x {\rm Pr}\{X\geq  x\}}\leq \left(1+\frac{1}{x^b}\right) \dfrac{d  }{d  x}\left[ \dfrac{1}{x f_X(x)}\right],
\end{align}
which means that the multiplicative inverse of the weighted (by $x$) right-tail probability, $x {\rm Pr}\{X\geq  x\}$,  of the  RV $X$ is  bounded   from above and bellow by  the derivative of the multiplicative inverse of its weighted  (by $x$)   PDF, $x f_X(x)$.

From \eqref{eq_5a3sdsd-2}, it is easy to see that the lower and upper bounds again converge to each other with rate $1/x^b$.

Lastly, note that the left-tail upper and lower bounds for a RV $X$ with support on $(-\infty,x_0])$ that satisfies assumptions 2)-6) when $x$ is negative, can be obtained straightforwardly from the derived right-tail upper and lower bounds simply by replacing $x$ with $-x$.

\section{Convergence Rates}

So far, I have presented upper and lower bounds on the right-tail, and corresponding conditions when these bounds hold. By checking whether the  conditions are met, one is able to see if the corresponding upper/lower bound holds or not. However, for a set of values $a$ in Theorem~\ref{thm_2} there is a corresponding set of upper bounds that jointly hold. Moreover, for a set of paired values $a$ and $b$ in Theorem~\ref{thm_4}, there is a corresponding set of lower bounds that jointly hold. Then, if one selects an upper bound from the first set and selects a lower bound from the second set, how close will these bounds be to the tail? This is the main topic of this  section. 

In the following, I refer to the process of selecting an upper bound from the first set and lower bound from the second set as upper and lower bound pairing, or simply pairing.
 Now,   pairing an upper bound with a corresponding lower bound can be done visually, by plotting them and then visually observing  the tightness between the two bounds. However, a more suitable method would be to utilize   some function that provides information about the tightness between the upper and lower bounds.  
 Note that  the tighter the upper and lower bounds are to each other, when they jointly hold, the tighter they are to the tail. This is made precise in the following.
 
 I propose to use the convergence rate between a potential upper bound and a potential lower bound  as a suitable function that provides information about the tightness between the two potential bounds. This is made precise in the following.

Let $P_U(x)$ and $P_L(x)$ be an upper and a lower bound of the right-tail ${\rm Pr}\{X\geq x\}$, respectively. Using $P_U(x)$ and $P_L(x)$, I define  the converge rate between the upper bound and the lower bound,   denoted by $R(x)$,  as
\begin{align}\label{eq_r1}
R(x)=\frac{P_U(x)}{P_L(x)} -1.
\end{align}
The function\footnote{Maybe a more accurate name for $R(x)$ would be divergence rate.} $R(x)$ provides information  about the speed of convergence between $P_U(x)$ and $P_L(x)$.
Note that the converge rate between the tail itself and its upper bound,  and the   converge rate between the tail itself and its   lower bound  can be bounded by $R(x)$ as
\begin{align}\label{eq_r2}
\max\left\{ \frac{P_U(x)}{{\rm Pr}\{X\geq x\}}-1, \frac{{\rm Pr}\{X\geq x\}}{P_L(x)}-1 \right\} \leq R(x),
\end{align}
where $R(x)$ is given by \eqref{eq_r1}. Hence, the   lower $R(x)$ is, the tighter the bounds are. In fact, $R(x)$ is zero if and only if $P_U(x)= P_L(x)= {\rm Pr}\{X\geq x\}$.

I first examine the convergence rate for RV $X$ with support on $([x_0,\infty)$, for $-\infty<x_0<\infty$. In that case, Theorems~\ref{thm_2} and \ref{thm_4} are applicable.
Now, the converge rate between the upper bound from Theorem~\ref{thm_2} for $a=a_U$, when it holds, and the lower bound from Theorem~\ref{thm_4} for $a=a_L$ and $b=b_L$, when it also holds, is given by
\begin{align}\label{eq_r3}
R(x)& = \dfrac{- \dfrac{y^{a_U} f_X^2(x)}{f_X(x)+y^{a_U} f'_X(x)}}{- \dfrac{y^{b_L}}{1+y^{b_L}}\dfrac{y^{a_L} f_X^2(x)}{f_X(x)+y^{a_L} f'_X(x)}} -1,  
\end{align}
where $y=x-x_0$. By simplifying \eqref{eq_r3}, I obtain
\begin{align}\label{eq_r4}
R(x)
&=\dfrac{1+y^{b_L}}{y^{a_L+b_L-a_U}}   \dfrac{ f_X(x)+y^{a_L} f'_X(x)} {f_X(x)+y^{a_U} f'_X(x)}   -1.
\end{align}
Now the function in \eqref{eq_r4} can be used for examining the tightness between a potential upper and a potential lower bound when $X$ has support on $([x_0,\infty)$, for $-\infty<x_0<\infty$. For example, if $a_U=a_L$   in \eqref{eq_r4}, then
\begin{align}\label{eq_r5}
R(x)   
 =  \dfrac{1}{y^{b_L}} =  \dfrac{1}{(x-x_0)^{b_L}}, 
\end{align}
and the smaller $R(x)$ is, the tighter the bounds are.

Now, if one seeks simplicity of the bounds, then a judicious choice would be pairing the upper bound from Corollary~\ref{cor_2}  with the lower bound from Corollary~\ref{cor_4} for a given $b=b_L$, if both bounds jointly hold. In that case, the rate of convergence is given by \eqref{eq_r5}, and the larger $b_L$ is, the tighter the bounds are. 

Another judicious choice for simplicity, in the case when the pairing between the bounds from  Corollaries~\ref{cor_2} and \ref{cor_4} is not possible since they do not hold jointly, is to pair  the bound from Corollary~\ref{cor_1}  with the bound from Corollary~\ref{cor_3} for a given $b=b_L$, if both bounds jointly hold. In that case, again the rate of convergence is given by \eqref{eq_r5}, and the larger $b_L$ is, the tighter the bounds are.

A third, maybe not so obvious choice for simplicity is pairing  the upper bound from Corollary~\ref{cor_1} with the lower bound from Corollary~\ref{cor_4} for a given $b=b_L$, if both bounds jointly hold.  In that case, $a_U=1$ and $a_L\to\infty$ holds, and thereby the convergence rate becomes
\begin{align}\label{eq_r6}
R(x)   
&=  \dfrac{1+y^{b_L}}{y^{b_L-1}} \frac{1}{y^{a_L}}   \dfrac{ f_X(x)+y^{a_L} f'_X(x)} {f_X(x)+y f'_X(x)}   -1,\nonumber\\
&\to \dfrac{1+y^{b_L}}{y^{b_L-1}} \frac{1}{y^{a_L}}   \dfrac{ y^{a_L} f'_X(x)} {f_X(x)+y f'_X(x)}   -1
 ,\;\; \textrm{as } a_L\to\infty,\nonumber\\
&= \dfrac{1+y^{b_L}}{y^{b_L-1}}    \dfrac{  f'_X(x)} {f_X(x)+y f'_X(x)}   -1
 ,\;\; \textrm{as } a_L\to\infty.
\end{align}
Depending on the RV, there might be cases when the convergence rate in \eqref{eq_r6} is better than that in \eqref{eq_r5}.

Now, the convergence rate when $X$ has support on $(-\infty,\infty)$  is straightforward. Specifically, in that case, the only possible choice is the upper bound from Corollary~\ref{cor_2} paired with the lower bound from Theorem~\ref{thm_5} for $b=b_L$, when they jointly hold. Thereby, the converge rate in that case is given by
\begin{align}\label{eq_r7}
R(x)& = \dfrac{- \dfrac{ f_X^2(x)}{  f'_X(x)}}{- \dfrac{(x-\mu)^{b_L}}{1+(x-\mu)^{b_L}} \dfrac{  f_X^2(x)}{  f'_X(x)}} -1 =  \dfrac{1+(x-\mu)^{b_L}} {(x-\mu)^{b_L}}-1 =  \dfrac{1 } {(x-\mu)^{b_L}}.
\end{align}
Hence, in this case, only the choice of $b_L$ matters. The larger $b_L$ is, the tighter the bounds are.

I use the converge rate for selecting the pairs of upper and lower bounds in the numerical examples in Section~\ref{sec-n}.

\section{Parameter Optimization}

In this section,  I propose methods for the optimizing  the parameter $a$ in the upper bounds and the parameter $b$ in the lower  bounds of the right-tail.

 \subsection{Optimizing Parameter $a$ In The Upper Bounds}

First note that the most general upper bound for a RV $X$ with support on $([x_0,\infty)$ is given in Theorem~\ref{thm_2}.
Now, note that the upper bound in Theorem~\ref{thm_2}, given by \eqref{eq_6a}, is a decreasing function of $a$ when $f'_X(x)<0$. Hence, for $f'_X(x)<0$, the following holds
\begin{align}\label{eq_oa1}
- \frac{ f_X^2(x)}{ f'_X(x)} \leq - \frac{(x-x_0)^a f_X^2(x)}{f_X(x)+(x-x_0)^a f'_X(x)}
\end{align}
when $x>x_0$.
This means that if the  conditions  in Corollary~\ref{cor_2}  are met, then there cannot be any further improvement of the upper bound  given in Theorem~\ref{thm_2} by optimizing $a$, simply because  in that case the upper bound in Corollary~\ref{cor_2} is the tightest possible upper bound. In other words, the optimal $a$ in that case is $a\to\infty$.

Now assume that the   conditions   in Corollary~\ref{cor_2}  are not met. Then, since the upper bound in Theorem~\ref{thm_2}, given by \eqref{eq_6a}, is a decreasing function of $a$, for $f'_X(x)<0$, the optimal $a$ is the largest possible $a$ for which conditions \eqref{eq_6b} and \eqref{eq_6c} are jointly satisfied. In practice, for a given $x$, this means that the optimal $a$, in this case, can be found by setting \eqref{eq_6c} as an equality,  solving it with respect to $a$, and choosing the largest possible solution for $a$ if there are multiple solutions.  Next, one needs to check if for this $a$ condition \eqref{eq_6b} still holds. If \eqref{eq_6b} does hold, then this is the optimal $a$. Otherwise, one should decrease $a$ or  increase $x$.   

\subsection{Optimizing Parameter $b$ In The Lower Bounds}

The easiest case to start  with is Theorem~\ref{thm_5}, which holds for RVs with support on $(-\infty, \infty)$. Now note that the lower bound in Theorem~\ref{thm_5}, given by \eqref{eq_4-1dfsd}, is an increasing function of $b$. Hence, the optimal $b$ in this case is the largest possible $b$ for which   conditions \eqref{eq_6c-nn3-n1-1-1-thm_7}  and \eqref{eq_6c-nn3-1-1-thm_7} jointly hold. As a result, for a given $x$ such that $f'_X(x)<0$ holds, the optimal $b$ in this case can be found by making \eqref{eq_6c-nn3-1-1-thm_7}  into an equality, solving it with respect to $b$, and then selecting the largest solution for $b$, if there  are multiple solutions.

For RVs with support on $([x_0, \infty)$, where $-\infty<x_0<\infty$, the optimal parameter $b$ can be found from the most general bound given by \eqref{eq_6a-nn1sd} in Theorem~\ref{thm_4}. Note that in \eqref{eq_6a-nn1sd} there are two parameters that can be optimized, $a$ and $b$, such that the bound is tightened when conditions \eqref{eq_6c-nn3-n1-1} and \eqref{eq_6c-nn3-1} are met. Since condition \eqref{eq_6c-nn3-n1-1} depends only on $a$, this condition is  met for any $a$ that satisfies 
\begin{align}\label{eq_aa1}
a >\frac{1}{x-x_0}\log_a\left(-\frac{f_X(x)}{f'_X(x)}\right)
\end{align}
when $f'_X(x)<0$ and $x>x_0$. Now there are two parameters $a$, that satisfies \eqref{eq_aa1}, and $b>0$ to make  \eqref{eq_6c-nn3-1} to be met with equality. Since this is an overparameterized system, i.e., one equation with two variables, I proposes the following sub-optimal approach. Set $a$ to be the optimal value that optimizes the upper bound in Theorem~\ref{thm_2}, for which I already explained above how it can be found. Then, optimize the parameter $b$ as follows. Since for fixed $x$ and $a$, the bound in Theorem~\ref{thm_4} is an increasing function of $b$,   find the optimal $b$ by setting \eqref{eq_6c-nn3-1} as an equality, solve it with respect to $b$, and choose the largest solution for $b$ in the case when there are multiple solutions.

\section{Numerical Examples}\label{sec-n}

In this section, I apply the derived upper and lower right-tail bounds to the following  RVs: the Gaussian, the non-central and central chi-square,  and the beta-prime, which have the following CDFs
\begin{align}
F_X(x)&=\frac{1}{2} \left(\mathrm{erf}\left(\frac{x-\mu }{\sqrt{2} \sigma }\right)+1\right),\label{eq_cdf1}\\
F_X(x)&=1-Q_{\frac{k}{2}}\left(\sqrt{\lambda },\sqrt{x}\right),\\
F_X(x)&=1-\frac{\Gamma \left(\frac{k}{2},\frac{x}{2}\right)}{\Gamma \left(\frac{k}{2}\right)},\\
F_X(x)&=\frac{B_{\frac{x}{x+1}}(\alpha ,\beta )}{B(\alpha ,\beta )},\label{eq_cdfk}
\end{align}
respectively, where $\mathrm{erf}(x)$, $Q_{M}(a,b)$, $\Gamma(a,b)$,  $\Gamma(a)$, $B_x(a,b)$, and $B(a,b)$ are the Gaussian error function, the Marcum-Q function, the gamma incomplete function, the  gamma function, the incomplete beta function, and the beta function, respectively. 

In the following numerical examples,   the parameters of the  CDFs in \eqref{eq_cdf1}-\eqref{eq_cdfk} are chosen almost at random in order to test the accuracy of the bounds.

\subsection{The Gaussian}

 Since the Gaussian distribution has unbounded support on $(-\infty,\infty)$, Corollary~\ref{cor_2} and  Theorem~\ref{thm_5} are applicable.

Let  the Gaussian distribution have mean $\mu=-1.7$ and standard deviation $\sigma=1.9$. 

\begin{figure}
\includegraphics[width=\textwidth]{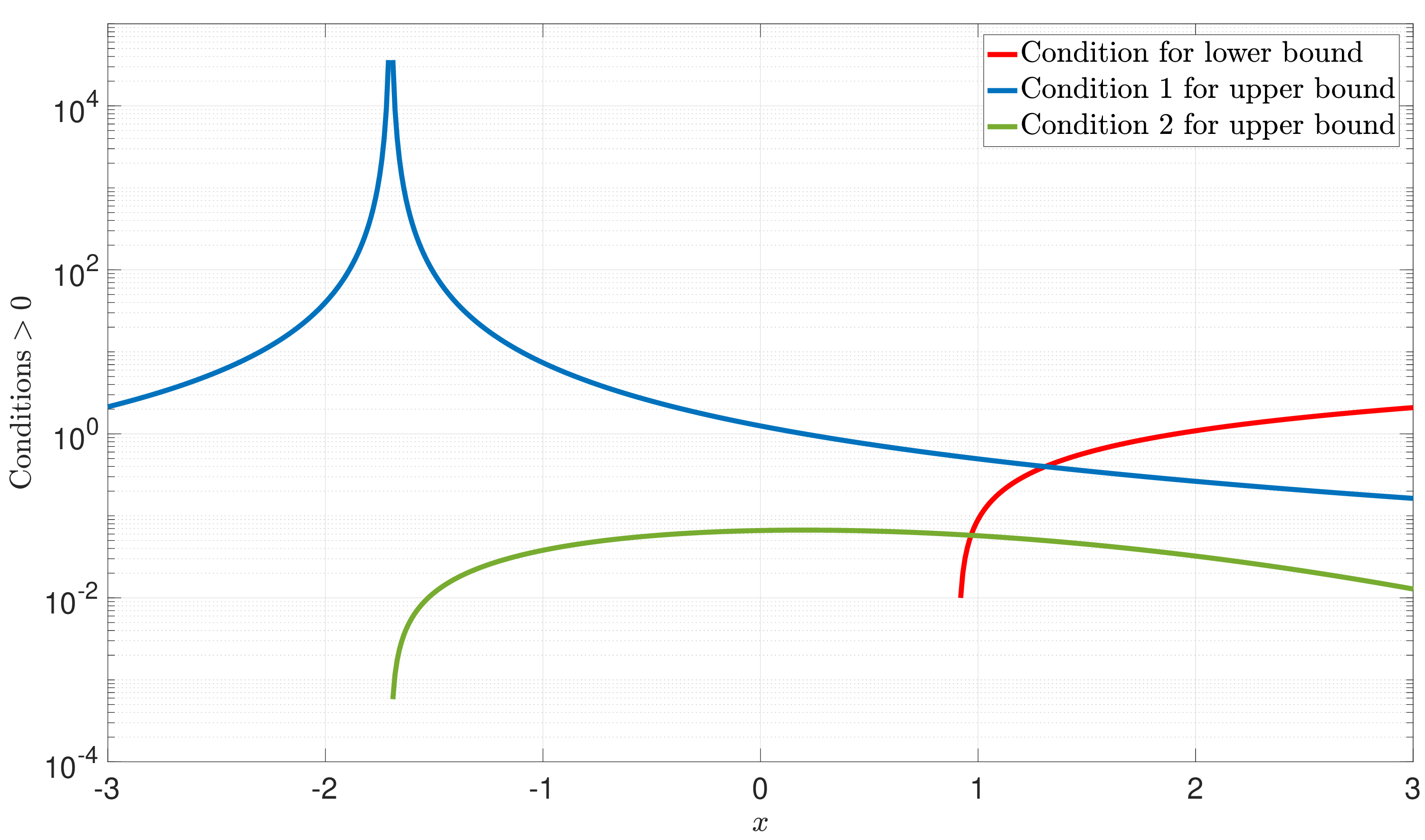}
\caption{The conditions.}
\label{fig_g1}
\end{figure}

\begin{figure}
\includegraphics[width=\textwidth]{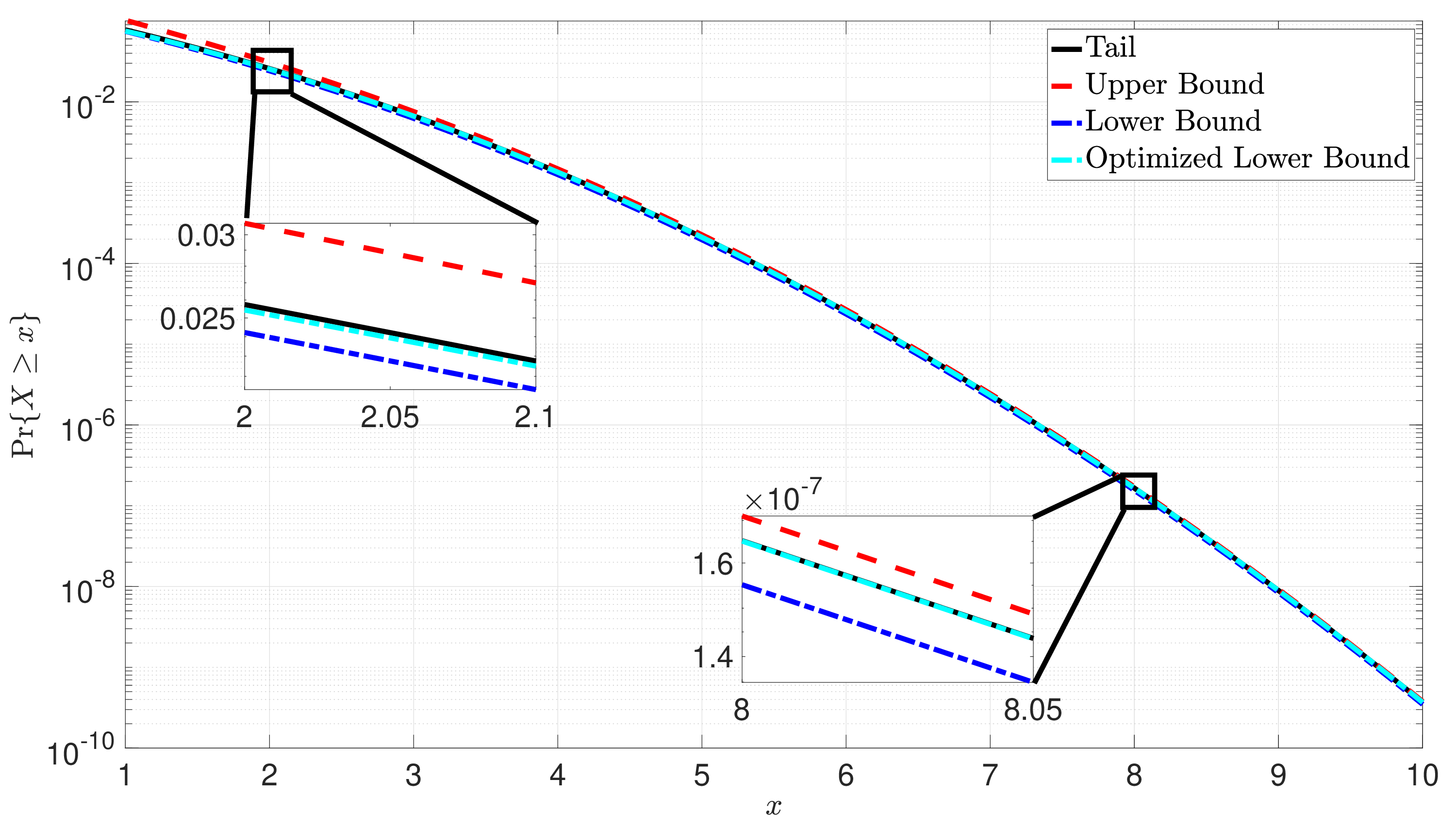}
\caption{The tail, the bounds, and the optimized lower bound for the Gaussia RV.}
\label{fig_g3}
\end{figure}

In Fig.~\ref{fig_g1}, I illustrate the necessary conditions for Corollary~\ref{cor_2} and  Theorem~\ref{thm_5} to hold jointly, when $b$ is not optimized. Specifically, in Fig.~\ref{fig_g1}, the green and blue lines are the negative values  of conditions \eqref{eq_5a-1-0} and \eqref{eq_5a-1}, which are the necessary conditions for the upper bound in Corollary~\ref{cor_2} to hold. Whereas, the red line is  condition \eqref{eq_6c-nn3-1-1-thm_7}, which is the main necessary condition for the lower bound in Theorem~\ref{thm_5} to hold. As this figure shows, it is necessary for $x$ to approximately be larger than one in order for both the lower and upper bounds to hold. As a result, the tail bounds are depicted for $x>1$ in Fig.~\ref{fig_g3}.

In Fig.~\ref{fig_g3}, I show the tail probability, the upper bound from Corollary~\ref{cor_2}, the lower bound from Theorem~\ref{thm_5} for $b=1$, and the lower bound from Theorem~\ref{thm_5} for optimized $b=1$. As can be seen from the figure, the tail, the upper bound, and the lower bounds are almost indistinguishable. Due to this, I have provided two zoomed regions, from where the bounds are more visible. As $x$ increases, the bounds become even tighter, which can be see from the second zoomed region. The optimized  lower bound is almost indistinguishable even in the first zoomed region, and almost identical to the tail bound in the second zoomed region.

If the expressions for the PDF and its derivative are plugged into the corresponding upper and lower bound for fixed $b$, the following bounds in closed form are obtained
\begin{align}\label{eq_g_a2}
\frac{\sigma  (x-\mu )^{b-1} }{\sqrt{2 \pi } \left((x-\mu )^b+1\right)} e^{-\frac{(x-\mu )^2}{2 \sigma ^2}} \leq {\rm Pr}\{X\geq x\}\leq \frac{\sigma  }{\sqrt{2 \pi } (x-\mu )} e^{-\frac{(x-\mu )^2}{2 \sigma ^2}},
\end{align}
which holds for $x>\mu$.
A good choice for $b$ is $b=1$.
Of course, the bounds in \eqref{eq_g_a2} hold when the corresponding conditions in Corollary~\ref{cor_2} and  Theorem~\ref{thm_5} jointly hold.  

As can be seen from \eqref{eq_g_a2}, the bounds obtained with this method are one of the tightest bounds for the Gaussian tail available in the literature, see \cite{1094433} for example. With this in mind,   the tightness of the bounds illustrated in Fig.~\ref{fig_g3} is not surprising.

\subsection{Non-Central and Central Chi-Square}

The non-central and central chi-squared RVs have been one of those RVs whose tails have been hard to bound. This is one of the reasons why these RVs have been selected for numerical evaluation in this section. But it turns out that these two RVs are also very interesting from the perspective of the proposed bounds, as will be discussed in the following.

Since the non-central and central chi-squared RVs have  support on $(0,\infty)$,   Theorem~\ref{thm_1} is applicable for the upper bound. Checking further if the optimal $a$ is $a\to\infty$, by checking if the conditions in Corollary~\ref{cor_2} hold, reveals that the bound in Corollary~\ref{cor_2} holds for $k\geq 2$, but it does not hold  for $k=1$. Instead,  for $k=1$, the conditions in Theorem~\ref{thm_1} hold. Therefore, the cases for $k=1$ and $k\geq 2$ have to be separated.

\subsubsection{The Case of $k\geq 2$}

   Since for $k\geq 2$, Corollary~\ref{cor_2} holds for the upper bound I adopt Corollary~\ref{cor_4}   for the lower bound.  Checking if the conditions in Corollary~\ref{cor_4} are met for $b=1$, it is revealed that they indeed are met. Hence, since I am seeking simplicity of the bounds, I adopt $b=1$.

\begin{figure}
\includegraphics[width=\textwidth]{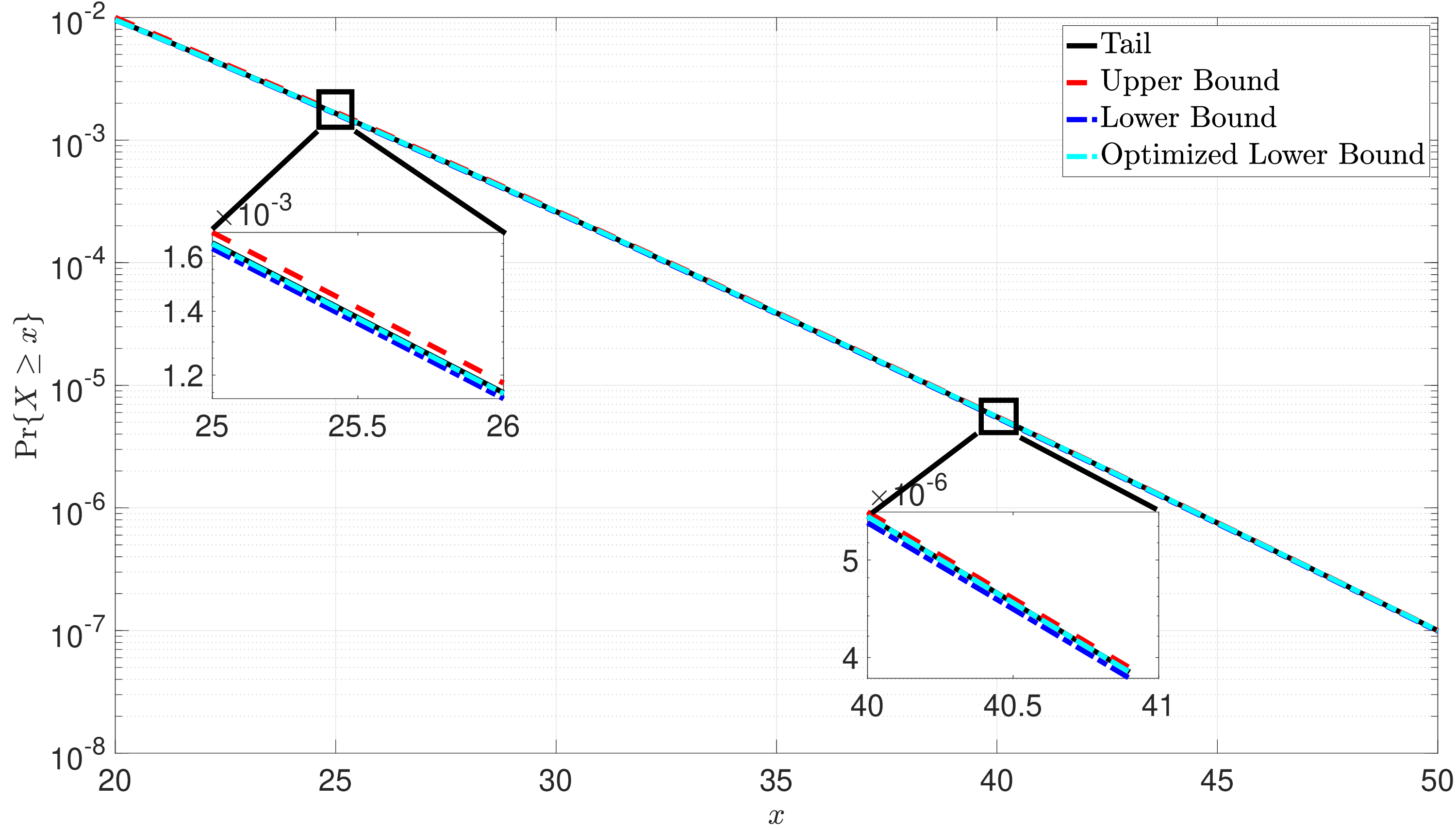}
\caption{The tail, the bounds, and the optimized lower bound for the non-central chi-squared RV.}
\label{fig_nccs}
\end{figure}

In Fig.~\ref{fig_nccs}, I plot the tail, the upper bound from Corollary~\ref{cor_2}, the lower bound from Corollary~\ref{cor_4} for $b=1$ and for optimized   $b$, for the case of  the non-central chi-squared RV with parameters  $k = 6$ and $\lambda = 0.2k$. Again, as can be seen from Fig.~\ref{fig_nccs}, the bounds are very tight.

If the expressions for the PDF and its derivative of the non-central chi-squared RV are plugged into the corresponding upper and lower bound for fixed $b$, the following bounds in closed form are obtained
\begin{align}\label{eq_ccsd-1}
&-\frac{x^b}{1+x^b}\frac{  \sqrt{\lambda  x} \left(\frac{x}{\lambda }\right)^{k/4} I_{\frac{k-2}{2}}\left(\sqrt{x \lambda
   }\right){}^2}{(k-x-2) I_{\frac{k-2}{2}}\left(\sqrt{x \lambda }\right)+\sqrt{\lambda  x} I_{\frac{k}{2}}\left(\sqrt{x \lambda }\right)} e^{-\frac{\lambda }{2}-\frac{x}{2}}
   \nonumber\\
    & \leq {\rm Pr}\{X\geq x\}
     \nonumber\\
    & \leq 
     -\frac{ \sqrt{\lambda  x} \left(\frac{x}{\lambda }\right)^{k/4} I_{\frac{k-2}{2}}\left(\sqrt{x \lambda }\right){}^2}{(k-x-2)
   I_{\frac{k-2}{2}}\left(\sqrt{x \lambda }\right)+\sqrt{\lambda  x} I_{\frac{k}{2}}\left(\sqrt{x \lambda }\right)}e^{-\frac{\lambda }{2}-\frac{x}{2}},
\end{align}
which hold only when $k\geq 2$ and when   $x$ in the denominator in \eqref{eq_ccsd-1} is large enough such that the denominator becomes negative, i.e., when
\begin{align}\label{eq_sdsdsgsg}
(k-x-2)
   I_{\frac{k-2}{2}}\left(\sqrt{x \lambda }\right)+\sqrt{\lambda  x} I_{\frac{k}{2}}\left(\sqrt{x \lambda }\right)<0.
\end{align}
 A good choice for $b$ in   \eqref{eq_ccsd-1}   is $b=1$.

Similar tightness of the bounds  as in Fig.~\ref{fig_nccs} are observed for the central chi-squared RV, for $k\geq 2$. Therefore, I omit this figure and   only provide the analytical bounds. Specifically,  if the expressions for the PDF and its derivative of the  central chi-squared RV are plugged into the corresponding upper and lower bound for fixed $b$, the following bounds in closed form are obtained
\begin{align}\label{eq_ccsd-2}
-\frac{2^{1-\frac{k}{2}} x^{b+\frac{k}{2}}}{\left(x^b+1\right) (k-x-2) \Gamma \left(\frac{k}{2}\right)} e^{-x/2} 
     \leq {\rm Pr}\{X\geq x\}
     \leq 
      -\frac{2^{1-\frac{k}{2}}  x^{k/2}}{(k-x-2) \Gamma \left(\frac{k}{2}\right)} e^{-x/2},
\end{align}
which hold only when $k\geq 2$ and when $x>k-2$. A good choice for $b$ in   \eqref{eq_ccsd-2}   is $b=1$.

\subsubsection{The Case of $k=1$}
The chi-squared distributions are defined for unit variance Gaussians. However, to generalize the result for $k=1$,  I use the general form of the Gaussian squared distribution, given by
\begin{align}\label{eq_eq_g_sq}
f_X(x)=\frac{1}{2} \frac{1}{\sqrt x}\frac{1}{ \sqrt{2 \pi\sigma ^2} } e^{-\frac{\left(\sqrt{x}+\mu \right)^2}{2 \sigma ^2}} +
\frac{1}{2} \frac{1}{\sqrt x}\frac{1}{ \sqrt{2 \pi\sigma ^2} } e^{-\frac{\left(\sqrt{x}-\mu \right)^2}{2 \sigma ^2}}.
\end{align}

In this case, Theorem~\ref{thm_1} is applicable for the upper bound. Since I am seeking for simplicity, I set $a=1$ in Theorem~\ref{thm_1}. I now have to look for the appropriate lower bound.  One possibility is the bound in Theorem~\ref{thm_3}, for $a=1$. However, if the conditions in  Theorem~\ref{thm_3} are checked, one would obtain that they hold for very small $b$, in which case the convergence rate between the upper and lower bounds would be very slow. Another possibility is the lower bound in Corollary~\ref{cor_4}. Checking the conditions in Corollary~\ref{cor_4}, it turns out that they hold for   $b=1$. Since I am seeking for simplicity of the bounds, the bound from Corollary~\ref{cor_4} is adopted for $b=1$.   

\begin{figure}
\includegraphics[width=\textwidth]{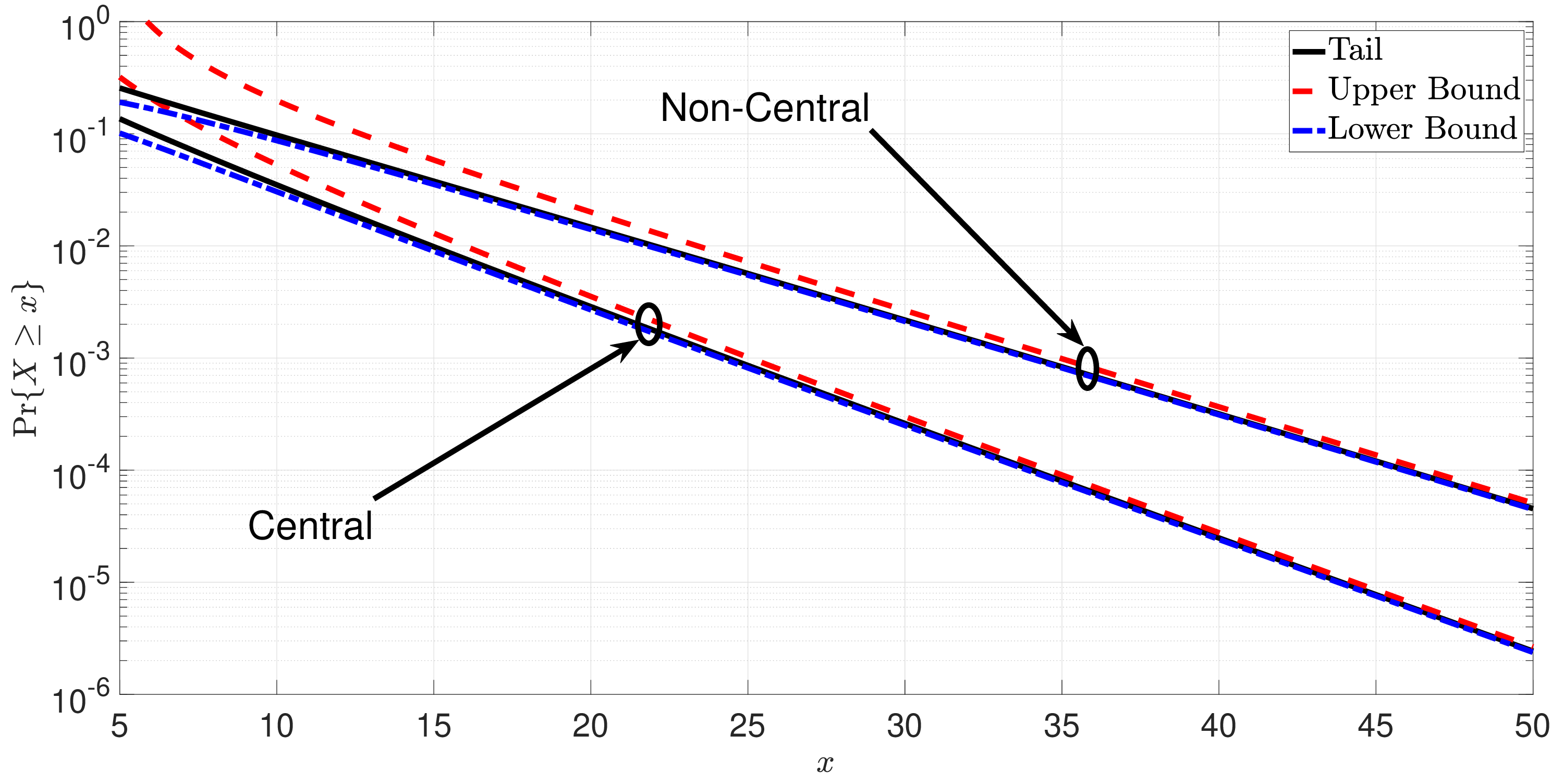}
\caption{The tail and the bounds  for the Gaussian squared RV.}
\label{fig_nccs2}
\end{figure}

In  Fig.~\ref{fig_nccs2}, I plot the tail, the upper bound from Theorem~\ref{thm_1} for $a=1$,   the lower bound  from Corollary~\ref{cor_4} for $b=1$,  for the Gaussian squared RV with 
 $\sigma=1.5$ and $\mu=0$ (central), and for $\mu=-1.2$ (non-central). As can be seen, from  Fig.~\ref{fig_nccs2}, the bounds are tight and becoming tighter as $x$ grows. The bounds using  optimized $a$ and $b$ are not shown as not to overcrowd the figure.

If the expressions for the PDF and its derivative of the   Gaussian squared RV are plugged into the corresponding upper and lower bound for fixed $a=1$ and $b=1$, the following bounds in closed form are obtained
\begin{align}\label{eq_ccsd-3}
&\frac{\sigma  x^{3/2} e^{-\frac{\left(\mu +\sqrt{x}\right)^2}{2 \sigma ^2}} \left(e^{\frac{2 \mu  \sqrt{x}}{\sigma ^2}}+1\right)^2}{\sqrt{2 \pi } (x+1) \left(x \left(e^{\frac{2 \mu  \sqrt{x}}{\sigma ^2}}+1\right)+\sigma ^2 \left(e^{\frac{2
   \mu  \sqrt{x}}{\sigma ^2}}+1\right)-\mu  \sqrt{x} \left(e^{\frac{2 \mu  \sqrt{x}}{\sigma ^2}}-1\right)\right)}
   \nonumber\\
    & \leq {\rm Pr}\{X\geq x\}\nonumber\\
    & \leq 
    \frac{\sigma  \sqrt{x} e^{-\frac{\left(\mu +\sqrt{x}\right)^2}{2 \sigma ^2}} \left(e^{\frac{2 \mu  \sqrt{x}}{\sigma ^2}}+1\right)^2}{\sqrt{2 \pi } \left(x \left(e^{\frac{2 \mu  \sqrt{x}}{\sigma ^2}}+1\right)- \sigma ^2 \left(e^{\frac{2
   \mu  \sqrt{x}}{\sigma ^2}}+1\right) -\mu  \sqrt{x} \left(e^{\frac{2 \mu  \sqrt{x}}{\sigma ^2}}-1\right)\right)},
\end{align}
which holds when  
 \begin{align}\label{eq_zdfcx}
 x \left(e^{\frac{2 \mu  \sqrt{x}}{\sigma ^2}}+1\right)- \sigma ^2 \left(e^{\frac{2
   \mu  \sqrt{x}}{\sigma ^2}}+1\right) -\mu  \sqrt{x} \left(e^{\frac{2 \mu  \sqrt{x}}{\sigma ^2}}-1\right) > 0,
 \end{align}
and is satisfied for large enough $x$.

\subsection{Beta-Prime}

The beta-prime RV is chosen since its PDF does not contain the exponential function.

Since the beta-prime RV has  support on $[0,\infty)$,   Theorem~\ref{thm_1} is applicable for the upper bound. Checking further if the optimal $a$ is $a\to\infty$, by checking if the conditions in Corollary~\ref{cor_2} hold, reveals that the bound in Corollary~\ref{cor_2}   does not hold. However, if   the conditions in Corollary~\ref{cor_1} are checked, then it is seen that the bound in Corollary~\ref{cor_1} holds. Therefore, for simplicity, I adopt the bound from Corollary~\ref{cor_1}. Now for the lower bound, checking  whether the conditions in  Corollary~\ref{cor_3} hold for $b=1$,  it turns out that they not hold. One needs to set $b$ to e.g. $b=4/5$ in order for the conditions in  Corollary~\ref{cor_3}, and the  corresponding lower bound, to hold,   for the chosen parameters. Another possibility is the lower bound in Corollary~\ref{cor_4}. Checking the conditions in Corollary~\ref{cor_4}, it turns out that they hold  even for $b\to\infty$, but even then the convergence rate between the lower and upper bounds is poorer than for the case when the lower bound from Corollary~\ref{cor_3} is chosen for $b=4/5$. Therefore, I pick $b=0.8$ in  Corollary~\ref{cor_3} as the lower bound.

\begin{figure}
\includegraphics[width=\textwidth]{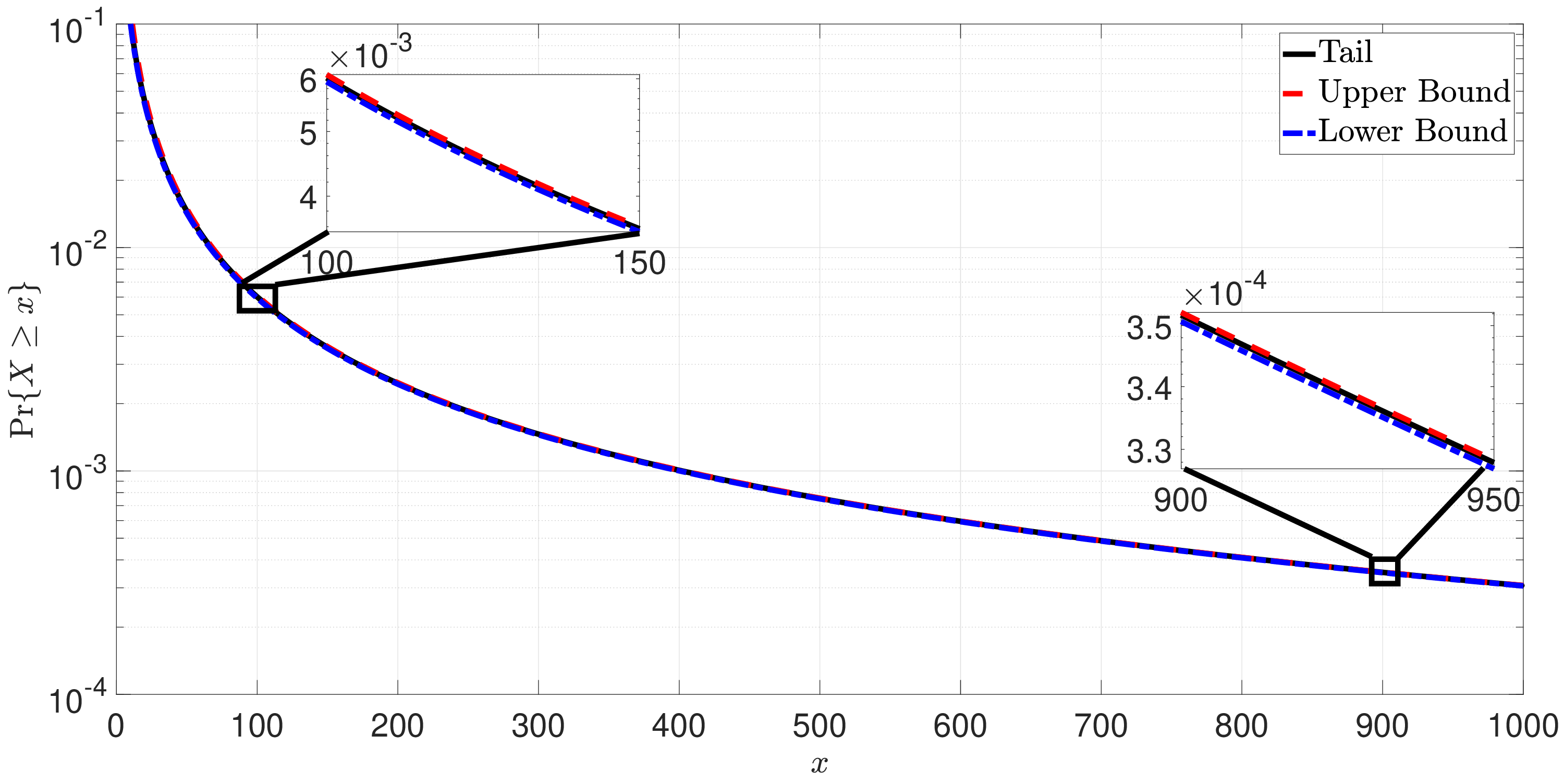}
\caption{The tail and the bounds for the beta-prime RV.}
\label{fig_nccs3}
\end{figure}

In Fig.~\ref{fig_nccs3}, I plot the tail, the upper bound from Corollary~\ref{cor_1},   the lower bound  from Corollary~\ref{cor_3} for $b=4/5$ for the beta-prime RV with parameters $\alpha=2.1$ and $\beta=1.3$. As can be seen, from  Fig.~\ref{fig_nccs2}, the bounds are again very tight.

If the expressions for the PDF and its derivative of the  beta-prime RV are plugged into the corresponding upper and lower bound for $b=4/5$, the following bounds in closed form are obtained
\begin{align}\label{eq_ccsd-4}
\frac{x^{\frac{4}{5}}}{ 1+ x^{\frac{4}{5}} }\frac{x^{\alpha} (x+1)^{-\alpha -\beta +1}}{ B(\alpha ,\beta ) (\beta  x-\alpha )}
     \leq {\rm Pr}\{X\geq x\}
     \leq 
      \frac{x^{\alpha } (x+1)^{-\alpha -\beta +1}}{B(\alpha ,\beta ) (\beta  x-\alpha )},
\end{align}
which holds  when  $x>\alpha/\beta$.

\newpage

 \appendix
 
 \subsection{Proof Of Theorem~\ref{thm_1}}\label{app_1}

First note, that proving an upper bound for the right-tail ${\rm Pr}\{X\geq x\}=1-F_X(x)$ is equivalent to proving a lower bound for $F_X(x)$. Specifically, if $F_X(x)\geq P_B(x)$ holds, then  the right-tail, ${\rm Pr}\{X\geq x\}$, is upper bounded as ${\rm Pr}\{X\geq x\}=1-F_X(x)\leq 1-P_B(x)$, where $P_B(x)$ is some function of $x$. In the following, I prove a lower bound for $F_X(x)$.

The integration by parts formula is given by 
\begin{align}\label{eq_6}
\int\limits_{0}^x u(t)  v'(t) dt  = v(t) u(t) \Bigg|_{0}^x- \int\limits_{0}^x  v(t) u'(t) dt.
\end{align}
On the other hand, $F_X(x)$ is given by\footnote{Note that expanding \eqref{eq_7} using Taylor series for some $x$, or expanding $f(t)$ for some $t$ and then integrating, does not lead to tight bounds.}
\begin{align}\label{eq_7}
F_X(x)=\int_{0}^x f_X(t) dt.
\end{align}
Applying the integration by parts on the integral in \eqref{eq_7}, where $u(t)=f_X(t)$ and $v'(t)=dt$, and  thereby  $u'(t)=f'_X(t)$ and $v(t)=t$, leads to
\begin{align}\label{eq_8}
F_X(x)=t f_X(t)\Bigg|_{0}^x- \int\limits_{0}^x  t f_X'(t) dt
\end{align}
or equivalently to
\begin{align}\label{eq_9}
F_X(x)=x f_X(x) - \int\limits_{0}^x  t f_X'(t) dt,
\end{align}
where $\lim\limits_{x\to 0} xf_X(x)=0$ is due to   assumption 6).
As a result of \eqref{eq_9}, the following holds
 \begin{align}\label{eq_10}
\dfrac{F_X(x)}{x f_X(x)- \int\limits_{0}^x  t f_X'(t) dt}=1.
\end{align}

 Now, taking the derivative with respect to $x$ on both sides of \eqref{eq_10}, I obtain
  \begin{align}\label{eq_11}
\dfrac{d }{d  x}\left[\dfrac{F_X(x)}{x f_X(x) - \int\limits_{0}^x  t f_X'(t)(t) dt}\right]=0.
\end{align}
On the other hand, the derivative in \eqref{eq_11} is obtained as
  \begin{align}\label{eq_12}
\dfrac{d }{d  x}\left[\dfrac{F_X(x)}{x f_X(x) - \int\limits_{0}^x  t f_X'(t) dt}\right]=\frac{F'_X(x)}{H_X(x)} -\dfrac{F_X(x)}{H_X^2(x)}{H'_X(x)},
\end{align}
where 
  \begin{align}
   F'_X(x)&=f_X(x),\label{eq_13a}\\
H_X(x)&= x f_X(x) - G_X(x),\label{eq_13b}\\
G_X(x)&= \int\limits_{0}^x  t f_X'(t) dt,\label{eq_13aa}\\
H'_X(x)&= f_X(x)+x f'_X(x) -   G'_X(x).\label{eq_13c} 
\end{align}
Note from \eqref{eq_13aa} that 
  \begin{align}
  G'_X(x)=x f'_X(x), \label{eq_13d}
\end{align}
  however, I keep $G'_X(x)$ as a separate auxiliary  function.

Inserting \eqref{eq_13a},  \eqref{eq_13b}, and \eqref{eq_13c} into \eqref{eq_12}, and then inserting \eqref{eq_12} into \eqref{eq_11}, I obtain
\begin{align}\label{eq_14}
\frac{f_X(x)}{x f_X(x) - G_X(x)} -\frac{F_X(x)}{(x f_X(x) - G_X(x))^2} \big(f_X(x)+x f'_X(x) -   G'_X(x)\big)=0.
\end{align}
From \eqref{eq_14}, I can obtain $F_X(x)$ as
\begin{align}\label{eq_15}
F_X(x)=\frac{ x f_X^2(x) }{f_X(x)+x f'_X(x)- G'_X(x)} - \frac{f_X(x) G_X(x)} {f_X(x)+x f'_X(x)- G'_X(x)}.
\end{align}
Note that if I insert \eqref{eq_13d} and \eqref{eq_13aa} into \eqref{eq_15} and simplify, I will obtain \eqref{eq_9}, which confirms the correctness of  \eqref{eq_15}.

I now insert  \eqref{eq_13d}  into \eqref{eq_15} but only into the second expression after the equality, and thereby obtain  
\begin{align}\label{eq_16}
F_X(x)&=\frac{x f_X^2(x) }{f_X(x)+x f'_X(x)- G'_X(x)} - \frac{f_X(x) G_X(x)} {f_X(x)+x f'_X(x)- x f'_X(x)}\nonumber\\
&=\frac{x f_X^2(x)    }{f_X(x)+x f'_X(x)- G'_X(x)} - \frac{f_X(x) G_X(x)} {f_X(x)},
\end{align}
which leads to the final expression of $F_X(x)$ given by
\begin{align}\label{eq_16a}
F_X(x)=\frac{x f_X^2(x)    }{f_X(x)+x f'_X(x)- G'_X(x)} -   G_X(x).
\end{align}

I now make my main claim in this proof, which is that, under certain assumptions that will be specified later on, the following holds
\begin{align}\label{eq_17}
\frac{x f_X^2(x) }{f_X(x)+x f'_X(x)- G'_X(x)} -   G_X(x) \geq \frac{x^a f_X^2(x)    }{f_X(x)+x^a f'_X(x)} +1,
\end{align}
when \eqref{eq_5} is met, where $a>0$.

Note that due to \eqref{eq_16a}, the claim that    \eqref{eq_17} holds is equivalent to claiming that the following holds under certain conditions 
\begin{align}\label{eq_18}
F_X(x)\geq \frac{x^a f_X^2(x)  }{f_X(x)+x^a f'_X(x)} +1.
\end{align}
Now, since the claim that \eqref{eq_17} holds, conditioned on \eqref{eq_5} being satisfied,  is equivalent to the claim that   \eqref{eq_18} holds, when \eqref{eq_5} is met, I prove  \eqref{eq_18}, which in fact is the main result of this theorem that I aim to prove.

I prove  \eqref{eq_18}, under certain conditions that will be specified at the end, by proving that the function
\begin{align}\label{eq_19}
D_X(x)=F_X(x)-  \frac{x^a f_X^2(x)  }{f_X(x)+x^a f'_X(x)} -1
\end{align}
is $i)$ a decreasing function, and $ii)$ that it  converges to zero as $x\to\infty$. The only possibility for $D_X(x)$ to be a decreasing function of $x$ and to $  D_X(x) \to 0$ as $x\to\infty$ is for $D_X(x)$ to be positive function of $x$ that converges to zero as $x\to\infty$.

To prove $i)$, I need to prove that, under the later specified conditions, the following holds 
\begin{align}\label{eq_20}
\dfrac{d  D_X(x) }{d  x}\leq  0.
\end{align}

I now prove $i)$. The derivative of $D_X(x)$ with respect to $x$ is given by
\begin{align}\label{eq_22}
\dfrac{d  D_X(x) }{d  x}&=
\frac{f_X(x)}{x \big(f_X(x)+x^a f'_X(x)\big)^2} \nonumber\\
&\times
\Big( (x-a x^a) f^2_X(x) 
- x^{2a+1} \big(f'_X(x)\big)^2
+x^{a+1}f_X(x)\big(f'_X(x)+x^a f''_X(x)\big) 
\Big)\leq 0.
\end{align}
Note that the term
\begin{align}\label{eq_23}
\frac{f_X(x)}{x \big(f_X(x)+x^a f'_X(x)\big)^2} 
\end{align}
in \eqref{eq_22} is positive. Hence, 
 \eqref{eq_20} holds
if and only if
\begin{align}\label{eq_25}
(x-a x^a) f^2_X(x) 
- x^{2a+1} \big(f'_X(x)\big)^2
+x^{a+1}f_X(x)\big(f'_X(x)+x^a f''_X(x)\big) \leq 0.
\end{align}
Hence, the necessary condition for $D_X(x)$ to be a decreasing function of $x$, i.e., \eqref{eq_20} to hold, is  \eqref{eq_25} to hold.
This concludes the proof of part $i)$. 

To prove $ii)$, I need to prove that, under   conditions  \eqref{eq_5}  and  \eqref{eq_25}, the following holds 
\begin{align}\label{eq_21}
\lim_{x\to\infty}    D_X(x)   = 0.
\end{align}
Taking the limit of $D_X(x)$, I obtain
\begin{align}
\lim_{x\to\infty}    D_X(x)       
&=    \lim_{x\to\infty}  \left(F_X(x)-  \frac{x^a f_X^2(x)  }{f_X(x)+x^a f'_X(x)} -1\right)  \label{eq_26b} \\
&=   \lim_{x\to\infty}  \big(F_X(x)-1\big )-  \lim_{x\to\infty}    \frac{x^a f_X^2(x)  }{f_X(x)+x^a f'_X(x)}   
\label{eq_26c} \\
&=   -  \lim_{x\to\infty}    \frac{x^a f_X^2(x)  }{f_X(x)+x^a f'_X(x)}    \label{eq_26d} \\
&=    -  \lim_{x\to\infty}    \frac{x^a f_X^2(x)  }{x^a f'_X(x)}  \label{eq_26e} \\
&=   -  \lim_{x\to\infty}    \frac{ f_X^2(x)  }{f'_X(x)}    
\label{eq_26f} \\
&     = 0.\label{eq_26h}
\end{align}
where  \eqref{eq_26b} is obtained by inserting \eqref{eq_19}, \eqref{eq_26e} holds  due to condition \eqref{eq_5} and \eqref{eq_4a}, and  \eqref{eq_26h} holds due to assumption 5). Note that I have not invoked condition \eqref{eq_25} for proving this limit. This concludes the proof of part $ii)$ part, and thereby concludes the proof of this theorem.

  \subsection{Proof Of Theorem~\ref{thm_2}}\label{app_2}
  
Let us define an auxiliary RV $Y$ as $Y=X-x_0$. Then the RV $Y$ has support on $([0,\infty)$ and as a result Theorem~\ref{thm_1} holds. Consequently, the following holds
   \begin{align} \label{eq_28}
{\rm Pr}\{Y\geq  y\}\leq - \frac{y^a f_Y^2(y)}{f_Y(x)+y^a f'_Y(y)}
\end{align}
under conditions 
\begin{align} 
&f_Y(y)+y^a f'_Y(y)<0\label{eq_5n-1}\\
&(y-a y^a) f^2_Y(y) 
- y^{2a+1} \big(f'_Y(y)\big)^2
+y^{a+1}f_Y(y)\big(f'_Y(y)+y^a f''_Y(y)\big) \leq 0.\label{eq_5an-2}
\end{align}

On the other hand,
 \begin{align} \label{eq_29}
{\rm Pr}\{X\geq  x\}={\rm Pr}\{Y+x_0\geq  x\} = {\rm Pr}\{Y\geq  x-x_0\} .
\end{align}
As a result of \eqref{eq_29} and \eqref{eq_28}, the following holds
 \begin{align} \label{eq_29a}
{\rm Pr}\{X\geq  x\}= {\rm Pr}\{Y\geq  x-x_0\} \leq - \frac{(x-x_0)^a f_Y^2(x-x_0)}{f_Y(x-x_0)+(x-x_0)^a f'_Y(x-x_0)}
\end{align}
under conditions 
\begin{align} 
&f_Y(x-x_0)+(x-x_o)^a f'_Y(x-x_0)<0,\label{eq_5n-1-a}\\
&
(x-x_0-a (x-x_0)^a) f^2_Y(x-x_0) 
- (x-x_0)^{2a+1} \big(f'_Y(x-x_0)\big)^2\nonumber\\
&
+(x-x_0)^{a+1}f_Y(x-x_0)\big(f'_Y(x-x_0)+(x-x_0)^a f''_Y(x-x_0)\big) \leq 0.\label{eq_5an-2-a}
\end{align}

Now, since
 \begin{align} \label{eq_29b}
{\rm Pr}\{X\leq  x\}={\rm Pr}\{Y+a\leq  x\} = {\rm Pr}\{Y\leq  x-a\} .
\end{align}
the following holds
 \begin{align} \label{eq_29c}
f_X(x)=f_Y(x-a).
\end{align}
Furthermore, due to \eqref{eq_29c}, the following holds
 \begin{align} \label{eq_29d}
f'_X(x)=f'_Y(x-a)
\end{align}
and
 \begin{align} \label{eq_29e}
f''_X(x)=f''_Y(x-a).
\end{align}

Inserting \eqref{eq_29c}, \eqref{eq_29d}, and \eqref{eq_29e} into \eqref{eq_29a}, \eqref{eq_5n-1-a}, and \eqref{eq_5an-2-a}, I obtain that
\begin{align} \label{eq_29a-1}
{\rm Pr}\{X\geq  x\} \leq - \frac{(x-x_0)^a f_X^2(x)}{f_X(x)+(x-x_0)^a f'_X(x)}
\end{align}
under conditions 
\begin{align} 
&f_X(x)+(x-x_o)^a f'_X(x)<0,\label{eq_5n-1-a-1}\\
&
(x-x_0-a (x-x_0)^a) f^2_X(x) 
- (x-x_0)^{2a+1} \big(f'_X(x)\big)^2\nonumber\\
&
+(x-x_0)^{a+1}f_X(x)\big(f'_X(x)+(x-x_0)^a f''_X(x)\big) \leq 0,\label{eq_5an-2-a-1}
\end{align}
which is the result I aimed to prove.

\subsection{Proof of Theorem~\ref{thm_3}}\label{app_3}
 
 The same method is followed as in the proof of Theorem~\ref{thm_1}, with the following difference. I now need to prove an upper bound for $F_X(x)$ in order to obtain a lower bound for the right-tail ${\rm Pr}\{X\geq x\}=1-F_X(x)$.

I prove an upper bound for $F_X(x)$, and thereby prove  \eqref{eq_6a-nn1}, under certain conditions that will be specified at the end, by proving that the function
\begin{align}\label{eq_19app-4}
D_X(x)=F_X(x)-\frac{x^b}{1+x^b}\frac{x^a f_X^2(x)}{f_X(x)+x^b f'_X(x)} -1
\end{align}
is $i)$ an increasing function and   $ii)$ that it  converges to zero as $x\to\infty$. The only possibility for $D_X(x)$ to be an increasing function of $x$ and to $  D_X(x) \to 0$ as $x\to\infty$ is for $D_X(x)$ to be negative function of $x$ that converges to zero as $x\to\infty$.

To prove $i)$, I need to prove that, under the later specified conditions, the following holds 
\begin{align}\label{eq_20app-4}
\dfrac{d  D_X(x) }{d  x}\geq  0.
\end{align}
 The derivative of $D_X(x)$ with respect to $x$ is given by
\begin{align}\label{eq_22app-4}
\dfrac{d  D_X(x) }{d  x}&=
\dfrac{f_X(x)}{ x\big(1+x^b\big)^2\big(f_X(x)+x^a f'_X(x)\big)^2}\nonumber\\
&
\times \Big( 
x\big(1+x^b\big)^2
- x^{a+b}(a+b+a x^b) f^2_X(x)
-x ^{2a+1}(x^{2b}-1) \big(f'_X(x)\big)^2\nonumber\\
&+ x^a f_X(x)\Big[ x(1+x^b)(2+x^b)f'_X(x)
+x^{a+b} \big(-b f'_X(x)+x(1+x^b)f''_X(x)\big)\Big]
 \Big)\geq 0.
\end{align}
Note that the term
\begin{align}\label{eq_23app-4}
\dfrac{f_X(x)}{ x\big(1+x^b\big)^2\big(f_X(x)+x^a f'_X(x)\big)^2}
\end{align}
in \eqref{eq_22app-4} is positive. Hence, \eqref{eq_20app-4} holds if and only if
\begin{align}\label{eq_25app-4}
&x\big(1+x^b\big)^2
- x^{a+b}(a+b+a x^b) f^2_X(x)
-x ^{2a+1}(x^{2b}-1) \big(f'_X(x)\big)^2\nonumber\\
&+ x^a f_X(x)\Big[ x(1+x^b)(2+x^b)f'_X(x)
+x^{a+b} \big(-b f'_X(x)+x(1+x^b)f''_X(x)\big)\Big]\geq 0
\end{align}
holds.

To prove $ii)$, I need to prove that, under the later specified conditions, the following holds 
\begin{align}\label{eq_21app-4}
\lim_{x\to\infty}  D_X(x)   = 0,
\end{align}
 which  straightforward, given the proof of Theorem~\ref{thm_1}, after noting that
\begin{align}
\lim_{x\to\infty} \frac{x^b}{1+x^b}=1.
\end{align}
for $b>0$.
 This   concludes the proof  of this theorem.

\subsection{Proof of Theorem~\ref{thm_5}}\label{app_4}
 
I follow the same approach as for the proof of Theorem~\ref{thm_3}. Thereby, 
\begin{align}
 D_X(x)=F_X(x)- \frac{(x-\mu)^b}{1+(x-\mu)^b} \frac{f_X^2(x)}{ f'_X(x)}-1,
\end{align}
and 
\begin{align}\label{eq_gg-1}
&\dfrac{d  D_X(x)}{d  x} 
= \frac{f_X(x)}{\big(1+(x-\mu)^b\big)^2}\nonumber\\
&
\times\Bigg(1-(x-\mu)^{2b}+
\frac{f_X(x)}{\big(f'_X(x)\big)^2}(x-\mu)^{b-1}
\Big(
-b f'_X(x)+\big(1+(x-\mu)^b\big)(x-\mu) f''_X(x)
\Bigg).
\end{align}
Noting that 
\begin{align}\label{eq_gg-2}
\frac{f_X(x)}{\big(1+(x-\mu)^b\big)^2}>0
\end{align}
for $x>\mu$, it follows that
\begin{align}\label{eq_gg-3}
&\dfrac{d  D_X(x)}{d  x} \geq 0
\end{align}
if an only if the following holds
\begin{align}\label{eq_gg-4}
1-(x-\mu)^{2b}+
\frac{f_X(x)}{\big(f'_X(x)\big)^2}(x-\mu)^{b-1}
\Big(
-b f'_X(x)+\big(1+(x-\mu)^b\big)(x-\mu) f''_X(x)
\geq 0.
\end{align}
Given the proof  of Theorem~\ref{thm_3}, the rest of the proof is straightforward and is therefore omitted.

\bibliographystyle{IEEEtran}
\bibliography{Citations}

\end{document}